\documentclass[12pt,a4]{article}

\usepackage{fixltx2e}
\usepackage{amsmath,amsthm,amssymb}
\usepackage[T1]{fontenc}
\usepackage[utf8]{inputenc}
\usepackage[hidelinks=true]{hyperref}

\theoremstyle{plain}

\makeatletter
\newtheorem*{rep@theorem}{\rep@title}
\newcommand{\newreptheorem}[2]{%
\newenvironment{rep#1}[1]{%
 \def\rep@title{#2~\ref*{##1}}%
 \begin{rep@theorem}}%
 {\end{rep@theorem}}}
\makeatother

\newtheorem{theorem}{Theorem}
\newreptheorem{theorem}{Theorem}

\newcommand\abs[1]{\lvert #1 \rvert}

\newcommand\parens[1]{\left(#1\right)}
\newcommand\phan{{\vphantom{1}}}
\newcommand\ZZ{\mathbb{Z}}
\newcommand\QQ{\mathbb{Q}}
\newcommand\class{\mathcal{K}}

\DeclareMathOperator{\Irr}{Irr}
\let\phi=\varphi

\newcommand\product[1]{\def\arg{#1}%
  \def\i{1}\arg
  \def\i{2}\arg
  \cdots
  \def\i{n}\arg
}

\newcommand\eqdots{
  \mathchoice
  {\mathrel{\setbox0=\hbox{$\displaystyle=$}
      \setbox1=\hbox{$\displaystyle \vdots$}
      \kern.5\wd0\kern-.5\wd1\copy1}}
  {\mathrel{\setbox0=\hbox{$\textstyle=$}
      \setbox1=\hbox{$\textstyle \vdots$}
      \kern.5\wd0\kern-.5\wd1\copy1}}
  {\mathrel{\setbox0=\hbox{$\scriptstyle=$}
      \setbox1=\hbox{$\scriptstyle \vdots$}
      \kern.5\wd0\kern-.5\wd1\copy1}}
  {\mathrel{\setbox0=\hbox{$\scriptscriptstyle=$}
      \setbox1=\hbox{$\scriptscriptstyle \vdots$}
      \kern.5\wd0\kern-.5\wd1\copy1}}
}

\allowdisplaybreaks[1]          

\title{Determination of Conjugacy Class Sizes from Products of Characters}
\author{Ivan Andrus\\
  \and P\'al Heged\H us}

\begin{document}

\maketitle

\begin{abstract}
  In~\cite{Robinson:2009}, Robinson showed that the character degrees are determined by knowing, for all $n$, the number of ways that the identity can be expressed as a product of $n$ commutators.
  Earlier, in~\cite{Strunkov:1991}, Strunkov showed that the existence of characters of $p$-defect $0$ can be determined by counting solutions to certain equations involving commutators and conjugates.

  In this paper, we prove analogs to Robinson's and Strunkov's theorems by switching conjugacy classes and characters.
  We show that counting the multiplicity of the trivial character in certain products of characters determines the conjugacy class sizes and existence of conjugacy classes with $p$-defect $0$.
\end{abstract}

\section{Introduction}

There are several known results relating the irreducible (complex) characters of $G$ to counting the number of ways an element can be expressed as a product of commutators.
Of interest in this paper are results of Robinson and Strunkov.
Robinson is able to determine all character degrees by knowing how many ways the identity can be represented as a product of commutators.
Strunkov, on the other hand, determines information about $p$-defects but only requires counting modulo $p$.

\begin{theorem}[\cite{Robinson:2009}]\label{thm:robinson}
  Given a finite group $G$, knowing the number of solutions, for $n=1,2,\dots,\abs{G}$, to equations $1=\product{[a_{\i},b_{\i}]}$ determines the character degrees (with multiplicity) of $G$.

  If the equations are instead of the form $1=\product{a_{\i}^{2}b_{\i}^{2}}$ then the degrees of real characters are obtained.
\end{theorem}

\begin{theorem}[\cite{Strunkov:1991}]\label{thm:strunkov}
  Let $G$ be a finite group and $S$ a Sylow $p$-subgroup of $G$.
  Let $f(x_{1},\dots,x_{k},u_{1},\dots,u_{l})$ be a function on $G$ that it is the product of elementary functions $[x_{i},x_{i+1}]$ and $u_{j}^{x_{j}}$ where $x_{i}\in G$ and $u_{j}\in S$.
  Moreover, let $k\geq2$ and let each variable appear in a single elementary function.
  Then $G$ has a $p$-block of defect $0$ if and only if the number of solutions to $g=f(x_{1},\dots,x_{k},u_{1},\dots,u_{l})$ is not divisible by $p\abs{S}^{l}$ for some $g\in G$.

  If $f$ contains at least one elementary multiplier of the form $x_{j}^{2}$ then  the existence of real characters with $p$-defect $0$ is determined instead.
\end{theorem}

In particular the existence of a $p$-block of defect $0$ is equivalent to the fact that the number of solutions to $g=[x_{1},x_{2}]$ is coprime to $p$ for some $g\in G$.

It is well known that there is a strong duality between conjugacy classes and irreducible characters.
Many theorems about conjugacy class sizes, for example, are similar to theorems about character degrees.
For instance, the prime power character degrees determine the order of the nilpotent residual~\cite{Cossey:1993}, and the prime power conjugacy class sizes determines the order of the hypercenter~\cite{Cossey:1992}.
In that vein, we prove analogous theorems to those of Robinson and Strunkov, showing that knowledge of character multiplication is enough to determine information about conjugacy class sizes.

\begin{theorem}\label{thm:class-sizes}
  Counting the multiplicity of the trivial character in all products of the form $\product{\chi_{\i}\overline{\chi_{\i}}}$ for $n=1,2,\dots,\abs{G}$ and $\chi_{i}\in\Irr(G)$, determines the conjugacy class sizes of $G$.

  If the products are of the form $\product{\chi_{\i}^{2}}$, then the sizes of real conjugacy classes are obtained.
\end{theorem}

Recall that a conjugacy class $\class$ has $p$-defect $0$ if $\abs{\class}_{p}=\abs{G}_{p}$.

\begin{theorem}\label{thm:defect-classes}
  Let $n\geq2$.
  A group $G$ has a conjugacy class of $p$-defect $0$ if and only
  \begin{equation*}
    \gamma_{n}(\phi) = \sum_{\chi_{i}\in\Irr(G)} [\phi,\product{\chi_{\i}\overline{\chi_{\i}}}]
  \end{equation*}
  is not divisible by $p$ for some irreducible character $\phi$.

  Likewise, $G$ has a real class of $p$-defect $0$ if and only if
  \begin{equation*}
    \delta_{n}(\phi) = \sum_{\chi_{i}\in\Irr(G)} [\phi,\product{\chi_{\i}^{2}}]
  \end{equation*}
  is coprime to $p$ for some $\phi\in\Irr(G)$.
\end{theorem}

\section{Proofs}

Our proof of Theorem~\ref{thm:class-sizes} is very different from that of Robinson, though the proof of Theorem~\ref{thm:defect-classes} is similar to Strunkov's.

\begin{reptheorem}{thm:class-sizes}
  Counting the multiplicity of the trivial character in all products of the form $\product{\chi_{\i}\overline{\chi_{\i}}}$ for $n=1,2,\dots,\abs{G}$ and $\chi_{i}\in\Irr(G)$, determines the conjugacy class sizes of $G$.

  If the products are of the form $\product{\chi_{\i}^{2}}$, then the sizes of real conjugacy classes are obtained.
\begin{proof}
  Let $\pi=\sum_{\chi\in\Irr(g)}\chi\overline{\chi}$ be the permutation character of $G$ acting on itself via conjugation, and let $\gamma_{n}(\phi)$ denote the multiplicity of $\phi$ in $\pi^n$.
  We begin by proving that the sequence $\{\gamma_i(1_G)\}$ determines the sizes of the centralizers, and therefore the conjugacy class sizes of $G$.

  First, note that $\pi(g)=\abs{C_{G}(g)}=\abs{G}/\abs{g^{G}}$.
  Consider the multiplicity of the trivial character in $\pi^{n}$:
  \begin{align*}
    [1_{G},\pi^{n}]
    & = \frac{1}{\abs{G}}\sum_{\class} \abs{\class} 1_{G}(g_{\class})\overline{\pi^{n}(g_{\class})} \\
    & = \sum_{\class} \frac{\abs{\class}}{\abs{G}} \parens{\frac{\abs{G}}{\abs{\class}}}^{n}        \\
    & = \sum_{\class} \parens{\frac{\abs{G}}{\abs{\class}}}^{n-1}
  \end{align*}
  where the sum is over all conjugacy classes $\class$, and $g_{\class}\in \class$.

  Let $a_{i}$ be the number of conjugacy classes of size $i$ and $C_{i}=\tfrac{\abs{G}}{i}$.
  Then the above equations can be reformulated
  \begin{align*}
    [1_{G},\pi]     & = \sum a_{i}                  \\
    [1_{G},\pi^{2}] & = \sum a_{i}C_{i}             \\
    [1_{G},\pi^{3}] & = \sum a_{i}^{\phan}C_{i}^{2} \\
                    & \eqdots
  \end{align*}
  Clearly $a_{i}=0$ for all $i>\abs{G}$, so the sums are finite.
  Suppose that the $[1_{G},\pi^{n}]$ are given and view each line as having variables $a_{i}$ with known coefficients $C_{i}^{j}$.
  Considering the first $\abs{G}$ equations gives a set of linear equations.
  There is a unique solution since the coefficient matrix is Vandermonde type, hence non-singular.
  Thus the sequence $[1_{G},\pi^{n}]$ for $n=1,2,\dots,\abs{G}$ determines the conjugacy class sizes of $G$.

  Now,
  \begin{align*}
    \pi^{n}
 & = \Big(\sum_{\chi\in\Irr(G)}\chi\overline{\chi}\Big)^{n}           \\
 & = \sum_{\chi_{i}\in\Irr(G)}\product{\chi_{\i}\overline{\chi_{\i}}} \\
 & = \sum_{\chi_{i}\in\Irr(G)}\abs{\product{\chi_{\i}}}^{2},
  \end{align*}
  thus the multiplicity $[1_{G},\pi^{n}]$ is equal to the multiplicity of $1_{G}$ in this last sum, as claimed.

  Now, let $\psi=\sum_{\chi\in\Irr(G)}\chi^2$.
  Note that
  \begin{align*}
    \psi(g) & =\sum_{\chi\in\Irr(G)}\chi(g)^2 \\
    & =\sum_{\chi\in\Irr(G)}\chi(g)\overline{\chi(g^{-1})} \\
    & =
    \begin{cases}
      \abs{C_{G}(g)} & \text{if $g$ is a real element;} \\
      0                  & \text{otherwise.}
    \end{cases}
  \end{align*}
  is essentially $\pi$ above, but restricted to real classes.

  The proof follows exactly the same as before except the sum
  \begin{equation*}
    [1_{G},\psi^{n}]
    = \sum_{\class=\class^{-1}} \parens{\frac{\abs{G}}{\abs{\class}}}^{n-1}
  \end{equation*}
  is over real classes, and therefore only determines their class sizes.
  Likewise,
  \begin{equation*}
    \psi^{n}
     = \Big(\sum_{\chi\in\Irr(G)}\chi^{2}\Big)^{n}
     = \sum_{\chi_{i}\in\Irr(G)}\parens{\product{\chi_{\i}}}^{2},
  \end{equation*}
  so $[1_{G},\psi^{n}]$ is the overall multiplicity of $1_{G}$ in all possible $n$-factor products of squares of irreducible characters, as claimed.
\end{proof}
\end{reptheorem}

Note that the multiplicity of an irreducible character in $\pi$ or $\psi$ is given by row sums in the character table
\begin{align*}
  [\phi,\pi]  & = \sum_{\class} \phi(g_{\class}); \\
  [\phi,\psi] & = \sum_{\class=\class^{-1}} \phi(g_{\class}).
\end{align*}

We now prove a partial analog of Strunkov's theorem.

\begin{reptheorem}{thm:defect-classes}
  Let $n\geq2$.
  A group $G$ has a class of $p$-defect $0$ if and only
  \begin{equation*}
    \gamma_{n}(\phi) = \sum_{\chi_{i}\in\Irr(G)} [\phi,\product{\chi_{\i}\overline{\chi_{\i}}}]
  \end{equation*}
  is not divisible by $p$ for some irreducible character $\phi$.

  Likewise, $G$ has a real class of $p$-defect $0$ if and only if
  \begin{equation*}
    \delta_{n}(\phi) = \sum_{\chi_{i}\in\Irr(G)} [\phi,\product{\chi_{\i}^{2}}]
  \end{equation*}
  is coprime to $p$ for some $\phi\in\Irr(G)$.
\begin{proof}
  Fix $n\geq2$. Consider the following equation for $\phi\in\Irr(G)$:
  \begin{align*}
    \gamma_{n}(\phi)
    & = [\phi,\pi^{n}]                                                                         \\
    & = \frac{1}{\abs{G}}\sum_{\class} \abs{\class} \phi(g_{\class})\overline{\pi^{n}(g_{\class})} \\
    & = \sum_{\class} \parens{\frac{\abs{G}}{\abs{\class}}}^{n-1}\phi(g_{\class}).
  \end{align*}
  It is a weighted sum of the row of $\phi$ in the character table.
  If $G$ has no class of $p$-defect $0$, then all coefficients $\abs{G}/\abs{\class}$ are divisible by $p$, and so are $\gamma_{n}(\phi)$ for every $\phi$.
  The weights at $p$-singular columns are $0$.

  Consider now the above equations mod $p$ and suppose that the sums are $0$ for all $\phi\in\Irr(G)$.
  That is, there is a mod $p$ combination of the $p$-regular columns of the character table that give the all zero column.
  Every Brauer character $\eta$ is a $\ZZ$-linear combination of regular characters~\cite[Corollary~2.16]{Navarro:1998}, and so the same linear combination of the columns of the Brauer character table is also zero.
  But the columns of the Brauer character table are linearly independent over $k$, an algebraically closed field of characteristic $p$~\cite[Theorem~1.19]{Navarro:1998}, and so a non-trivial combination of them cannot be zero.
  Consequently all weights are divisible by $p$, as required.

  All the same arguments hold for real classes when $\pi$ is replaced by $\psi$.
\end{proof}
\end{reptheorem}

Since $\pi^{n}\psi^{m}=\psi^{n+m}$ information about real classes (for both theorems) is determined by counting over products of the form
\begin{equation*}
  \product{\abs{\chi_{\i}}^{2}}\chi_{n+1}^{2}\cdots\chi_{n+m}^{2}
\end{equation*}
as long as $m\geq1$ and $n+m\geq2$.

\section{Remarks}

Theorem~\ref{thm:defect-classes} is not a complete analog of Strunkov's because we only include an analog for commutators (or squares), not for conjugates of elements of $S$.
One must first determine what the analog of elements of $S$ would be.
The conjugates of elements in $S$ comprise the conjugacy classes of $p$-elements.

Let $\varepsilon$ be a primitive $|G|$-th root of unity and $R$ the algebraic integers of $\QQ[\varepsilon]$.
Let $M$ be a maximal ideal of $R$ containing the prime $p$.
The element $g\in G$ is a $p$-element if and only if, for all $\chi\in\Irr(G)$, $\chi(g)\equiv\chi(1)\mod M$.
See~\cite[8.20]{Isaacs:1976}.

An analog for characters based on similar congruence criteria is that $\chi$ is in the principal block $B_0$ of $G$ if and only if $\frac{\chi(g)\abs{g^{G}}}{\chi(1)}\equiv\abs{g^{G}}\mod M$ for all $g\in G$.
The analog of $u_{j}^{x_{j}}$ would be $|\chi|^2\phi$.
However, the corresponding result does not hold.

As a counter example, consider $S_{3}$ for $p=3$.
It has a unique $3$-block $B_{0}$, and the class of transpositions is of $3$-defect $0$.
Calculating
\begin{equation*}
  \gamma(\psi)
   = \sum_{\stackrel{\chi_{1},\chi_{2},\chi_{3}\in\Irr(G)}{\scriptscriptstyle\phi\in\Irr(B_{0})}}
  [\psi, \abs{\chi_{1}\chi_{2}}^{2}\abs{\chi_{3}}^{2}\phi]
\end{equation*}
(the analog of $[x_{1},x_{2}]u^{x_{3}}$) gives results divisible by $9=3|G|_3$ for all $\psi\in\Irr(G)$.
One might consider using $\sum_{\phi\in\Irr(B_{0})}\phi(1)^{2}$ (or it's $p$-part) in place of $\abs{S}=\abs{G}_{p}$, but this causes other groups to not satisfy the theorem (e.g., the dihedral group of order $12$).
Thus it seems that, if it exists, the correct analog is not the principal block.

\bibliographystyle{alpha}
\bibliography{class-sizes}

\begin{thebibliography}{CHM92}

\bibitem[CH93]{Cossey:1993}
John Cossey and Trevor Hawkes.
\newblock Computing the order of the nilpotent residual of a finite group from
  knowledge of its group algebra.
\newblock {\em Arch. Math. (Basel)}, 60(2):115--120, 1993.

\bibitem[CHM92]{Cossey:1992}
John Cossey, Trevor Hawkes, and Avinoam Mann.
\newblock A criterion for a group to be nilpotent.
\newblock {\em Bulletin of the London Mathematical Society}, 24(3):267--270,
  1992.

\bibitem[Isa76]{Isaacs:1976}
Irving~Martin Isaacs.
\newblock {\em Character theory of finite groups}.
\newblock Academic Press, New York, 1976.

\bibitem[Nav98]{Navarro:1998}
Gabriel Navarro.
\newblock {\em Characters and Blocks of Finite Groups}.
\newblock Number 250 in London Mathematical Society Lecture Notes Series.
  Cambridge University Press, 1998.

\bibitem[Rob09]{Robinson:2009}
Geoffrey~R. Robinson.
\newblock Characters and the commutator map.
\newblock {\em Journal of Algebra}, 321(11):3521 -- 3526, 2009.
\newblock Special Issue in Honor of Gus Lehrer.

\bibitem[Str91]{Strunkov:1991}
S.~P. Strunkov.
\newblock Existence and the number of {$p$}-blocks of defect {$0$} in finite
  groups.
\newblock {\em Algebra i Logika}, 30(3):355--368, 381, 1991.

\end{thebibliography}

\end{document}